\documentclass[12pt,draft]{amsart}
\def\NZQ{\Bbb}               
\def\NN{{\NZQ N}}

\def\AA{{\NZQ A}}

\def\opn#1#2{\def#1{\operatorname{#2}}} 
\opn\chara{char}
\opn\length{\ell}
\opn\pd{pd}
\opn\rk{rk}
\opn\projdim{proj\,dim}
\opn\injdim{inj\,dim}
\opn\rank{rank}
\opn\depth{depth}
\opn\grade{grade}
\opn\height{height}
\opn\embdim{emb\,dim}
\opn\codim{codim}

\opn\Tr{Tr}
\opn\bigrank{big\,rank}
\opn\superheight{superheight}\opn\lcm{lcm}
\opn\trdeg{tr\,deg}%
\opn\reg{reg}
\opn\lreg{lreg}
\opn\skel{skel}
\opn\com{com}
\opn\div{div}
\opn\Div{Div}
\opn\cl{cl}
\opn\Cl{Cl}
\opn\Spec{Spec}
\opn\Supp{Supp}
\opn\supp{supp}
\opn\Sing{Sing}
\opn\Ass{Ass}
\opn\Ann{Ann}
\opn\Rad{Rad}
\opn\Soc{Soc}
\opn\Ker{Ker}
\opn\Coker{Coker}
\opn\Im{Im}
\opn\Hom{Hom}
\opn\Tor{Tor}
\opn\Ext{Ext}
\opn\End{End}
\opn\Aut{Aut}
\opn\id{id}

\opn\nat{nat}
\opn\pff{pf}
\opn\Pf{Pf}
\opn\GL{GL}
\opn\SL{SL}
\opn\mod{mod}
\opn\ord{ord}
\opn\aff{aff}
\opn\con{conv}
\opn\relint{relint}
\opn\st{st}
\opn\lk{lk}
\opn\cn{cn}
\opn\core{core}
\opn\vol{vol}
\opn\link{link}
\opn\gr{gr}

\def\pot#1#2{#1[\kern-0.28ex[#2]\kern-0.28ex]}

\opn\dirlim{\underrightarrow{\lim}}
\opn\inivlim{\underleftarrow{\lim}}
\let\to=\rightarrow

\def\Implies{\ifmmode\Longrightarrow \else
     \unskip ${}\Longrightarrow{} $ \ignorespaces\fi}
\def\implies{\ifmmode\Rightarrow \else
     \unskip ${}\Rightarrow{} $ \ignorespaces\fi}
\def\iff{\ifmmode\Longleftrightarrow \else
     \unskip ${}\Longleftrightarrow{} $  \ignorespaces\fi}

\let\:=\colon

\newtheorem{Theorem}{Theorem}[section]

\newtheorem{Corollary}[Theorem]{Corollary}
\newtheorem{Proposition}[Theorem]{Proposition}
\newtheorem{Example}[Theorem]{Example}

\theoremstyle{definition}

\let\epsilon\varepsilon
\let\phi=\varphi
\let\kappa=\varkappa
\def\AA{{\mathcal A}}

\opn\initial{in}
\opn\inim{inm}
\opn\rev{rev}
\opn\Gin{Gin}
\opn\Lex{Lex}
\opn\Shift{Shift}
\opn\shift{shift}
\opn\rate{rate}
\opn\Mon{Mon}
\opn\lex{lex}
\opn\rev{rev}
\opn\red{red}
\opn\max{max}
\opn\min{min}
\opn\initial{in}
\opn\Ker{Ker}
\opn\GL{GL}
\opn\proj{proj}
\begin{document}
\title{
The depth of an ideal with a given Hilbert function
}
\author{Satoshi Murai}
\address{
Department of Pure and Applied Mathematics\\
Graduate School of Information Science and Technology\\
Osaka University\\
Toyonaka, Osaka, 560-0043, Japan\\
}
\email{s-murai@ist.osaka-u.ac.jp}
\thanks{The first author is supported by JSPS Research Fellowships for Young Scientists}

\author{Takayuki Hibi}
\address{
Department of Pure and Applied Mathematics\\
Graduate School of Information Science and Technology\\
Osaka University\\
Toyonaka, Osaka, 560-0043, Japan\\
}
\email{hibi@math.sci.osaka-u.ac.jp
}


\maketitle

\begin{abstract}
Let $A = K[x_1, \ldots, x_n]$ denote the polynomial ring
in $n$ variables over a field $K$ with each $\deg x_i = 1$.
Let $I$ be a homogeneous ideal of $A$ with $I \neq A$
and $H_{A/I}$ the Hilbert function of 
the quotient algebra $A / I$. 
Given a numerical function 
$H : \NN \to \NN$ satisfying
$H=H_{A/I}$ for some homogeneous ideal $I$ of $A$,
we write $\mathcal{A} _H$
for the set of those integers $0 \leq r \leq n$ 
such that there exists a homogeneous ideal $I$ of $A$
with $H_{A/I} = H$ and with
$\depth A / I = r$.
It will be proved that
one has either $\mathcal{A}_H = \{ 0, 1, \ldots, b \}$
for some $0 \leq b \leq n$
or $|\AA_H| = 1$.
\end{abstract}
\section*{Introduction}
Let $A = K[x_1, \ldots, x_n]$ denote the polynomial ring
in $n$ variables over a field $K$ with each $\deg x_i = 1$.
Let $I$ be a homogeneous ideal of $A$ with $I \neq A$
and $H_R$ the Hilbert 
function of the quotient algebra $R = A / I$.
Thus $H_R(q)$, $q = 0, 1, 2, \ldots$,
is the dimension of the subspace of $R$ 
spanned over $K$ by the homogeneous elements of $R$ 
of degree $q$. 
A classical result \cite[Theorem 4.2.10]{BruHer}
due to Macaulay
guarantees that, given a numerical function 
$H : \NN \to \NN$,
where $\NN$ is the set of nonnegative integers, 
there exists a homogeneous ideal $I$ of $A$ with
$I \neq A$ such that $H$ is the Hilbert function
of the quotient algebra $R = A / I$ 
if and only if $H(0) = 1$, $H(1) \leq n$
and $H(q+1) \leq H(q)^{\langle q \rangle}$ 
for $q = 1, 2, \ldots$, where
$H(q)^{\langle q \rangle}$ is defined in
\cite[p.\ 161]{BruHer}.

Given a numerical function 
$H : \NN \to \NN$ satisfying $H(0) = 1$, $H(1) \leq n$
and $H(q+1) \leq H(q)^{\langle q \rangle}$ 
for $q = 1, 2, \ldots$,
we write $\AA_H$ 
for the set of those integers $0 \leq r \leq n$ 
such that there exists a homogeneous ideal $I$ of $A$
with $H_{A/I} = H$ and with
$\depth A / I = r$.
We will show that,
given a numerical function 
$H : \NN \to \NN$ satisfying $H(0) = 1$, $H(1) \leq n$
and $H(q+1) \leq H(q)^{\langle q \rangle}$ 
for $q = 1, 2, \ldots$,
one has (i) $\AA_H = \{ n - \delta \}$ if
$H$ is of the form (\ref{Boston}) of
Proposition \ref{hilbert}
and (ii)
$\AA_H = \{ 0, 1, \ldots, b \}$ for some $b \geq 0$
if $H$ cannot be of the form (\ref{Boston}).
The statement (i) will be proved in Theorems \ref{critical},
and the statement (ii) will be proved in Theorem \ref{noncritical}.
Also, we will introduce a way to determine the integer $b=\max\AA_H$ from $H$
in Theorem \ref{added}.

\section{Universal lexsegment ideals}
Let $A = K[x_1, \ldots, x_n]$ denote the polynomial ring
in $n$ variables over a field $K$ with each $\deg x_i = 1$
and $A_{[m]}=K[x_1,\dots,x_{n+m}]$, where $m$ is a positive integer.
Work with the lexicographic order $<_{\lex}$ on $A$
induced by the ordering $x_1 > x_2 > \cdots > x_n$
of the variables. 
Write, as usual, $G(I)$ for the (unique) minimal system 
of monomial generators of a monomial $I$ of $A$.
Recall that a monomial ideal $I$ of $A$ is 
{\em lexsegment} 
if, for a monomial $u$ of $A$ belonging to $I$ 
and for a monomial $v$ of $A$ with $\deg u = \deg v$
and with $v >_{\lex} u$, one has $v \in I$.
A lexsegment ideal $I$ of $A$ is called {\em universal
lexsegment} (\cite{BNT})
if, for any integer $m \geq  1$,
the monomial ideal $I A_{[m]}$ of the polynomial ring $A_{[m]}$
is lexsegment.
In other words, a universal lexsegment ideal of $A$ is 
a lexsegment ideal $I=(u_1,\dots,u_t)$ of $A$ which remains being 
lexsegment if we regard  $I=(u_1,\dots,u_t)$ as an ideal of the polynomial ring
$A_{[m]}$ for  all $m \geq 1$.

\begin{Example}
{\em
(a) The lexsegment ideal
$(x_1^2, x_1 x_2^2)$ of $K[x_1, x_2]$ 
is universal lexsegment.
In fact, 
the ideal $(x_1^2, x_1 x_2^2)$ of 
$K[x_1, \ldots, x_m]$ is lexsegment
for all $m \geq  2$.

(b) The lexsegment ideal
$(x_1^3, x_1^2 x_2, x_1 x_2^2)$ of $K[x_1, x_2]$ 
cannot be universal lexsegment.
Indeed, since
$x_1 x_2^2 <_{\lex} x_1^2 x_3$
in $K[x_1, x_2, x_3]$,
the ideal $(x_1^3, x_1^2 x_2, x_1 x_2^2)$ of 
$K[x_1, x_2, x_3]$ is not lexsegment.
}
\end{Example}

\begin{Proposition}
\label{lexsegment}\
\begin{itemize}
\item[(a)]
Let $I$ be a lexsegment ideal of $A$ with
$G(I) = \{ u_1, \ldots ,u_\delta\}$
where $\deg u_1 \leq  \cdots \leq \deg u_\delta$
and where 
$u_{i+1} <_{\lex} u_i$
if $\deg u_i = \deg u_{i+1}$.
Let $s_1 = \deg u_1 - 1$ and
$s_i = \deg u_i - \deg u_{i-1}$
for $i = 2, 3, \ldots,\delta$.
Then, for $k \leq n$, one has
\[
u_k = x_1^{s_1} x_2^{s_2} \cdots x_k^{s_k + 1}.
\]
\item[(b)]
Given an integer $1 \leq \delta \leq n$ together with
a sequence of integers $1 \leq e_1 \leq \cdots \leq e_\delta$,
there is a lexsegment ideal $I$ of $A$
with $G(I) = \{ u_1, \ldots, u_\delta \}$ such that
$\deg u_i = e_i$ for $i = 1, \ldots, \delta$.
\end{itemize}
\end{Proposition}

\begin{proof}
(a)
Since $u_1 = x_1^{\deg u_1}$,
one has $u_1 = x_1^{s_1 + 1}$.
Let $1 < k \leq \min\{n,\delta\}$ and suppose that 
$u_{k-1} 
= x_1^{s_1} x_2^{s_2} \cdots x_{k-1}^{s_{k-1} + 1}$.
Since the ordering of $u_1,u_2,\dots,u_\delta$ implies that
the monomial ideal $(u_1, \ldots, u_{k-1})$ is lexsegment,
the smallest monomial with respect to $<_{\lex}$ of degree $\deg u_k$
belonging to $(u_1, \ldots, u_{k-1})$
is $u_{k-1} x_n^{s_k}$.
Since $u_k$ is the largest monomial with respect to $<_{\lex}$ which satisfies
$\deg u_{k} = \deg (u_{k-1} x_n^{s_k})$ and $u_k <_{\lex} u_{k-1}x_n^{s_k}$,
we have $u_k = (u_{k-1} / x_{k-1}) x_k^{s_k+1}$.
Thus $u_k =
x_1^{s_1} x_2^{s_2} \cdots x_{k-1}^{s_{k-1}} x_k^{s_k+1}$,
as desired.

(b)
This can be easily done by induction on $\delta$.
Let $\delta \leq n$ and suppose that $J$ is a lexsegment
ideal of $A$ 
with $G(J) = \{ u_1, \ldots, u_{\delta - 1} \}$
such that 
$\deg u_i = e_i$ for $i = 1,2, \ldots, \delta - 1$.
Then by (a) 
we have $G(J) \subset K[x_1, \ldots, x_{\delta-1}]$.  
Hence $x_{\delta}^{e_\delta} \not\in J$.
Thus there exists a monomial of degree $e_\delta$ which does not belong to $J$.
Let $u_{\delta}$ be the largest monomial of degree $e_{\delta}$
with respect to $<_{\lex}$ which does not belong to $J$.
Then $(u_1, \ldots, u_{\delta - 1}, u_\delta)$
is a lexsegment ideal of $A$.
\end{proof}

\begin{Corollary}   
\label{universal}
A lexsegment ideal $I$ of $A$ is universal lexsegment
if and only if $|G(I)| \leq n$, where
$|G(I)|$ is the number of monomials belonging to $G(I)$. 
\end{Corollary}

\begin{proof}
Let $G(I) = \{ u_1,\ldots,u_\delta \}$,
where $\deg u_1 \leq  \cdots \leq \deg u_\delta$.
It $\delta \geq n+1$, then $I A_{[1]}$ is not a lexsegment ideal of $A_{[1]}$
since Proposition \ref{lexsegment} (a)
says that, for any lexsegment ideal $J$ of $A_{[1]}$ with $|G(J)| \geq n+1$,
there exists a generator $v \in G(J)$ such that $x_{n+1}$ divides $v$.
Thus $I$ is not a universal lexsegment if $\delta \geq n+1$.

Assume that $\delta \leq n$.
For any positive integer $m$,
Proposition \ref{lexsegment} (b) says that there exists the lexsegment ideal $J$ of $A_{[m]}$
such that $G(J)= \{ v_1,\ldots,v_\delta \}$ satisfies $\deg v_i =\deg u_i$ for $i=1,2,\dots,\delta$.
Then Proposition \ref{lexsegment} (a) says that $G(I)=G(J)$.
Thus $I A_{[m]}$ is a lexsegment ideal of $A_{[m]}$ for all $m \geq 1$ if $\delta \leq n$.
\end{proof}

For any monomial $u$ of $A$,
let $m(u)$ be the biggest integer $1 \leq i \leq n$ 
for which $x_i$ divides $u$.
A monomial ideal $I$ of $A$ is said to be \textit{stable}
if $u \in I$ implies $(x_{q}/x_{m(u)}) u \in I$ for any $1 \leq q <m(u)$.
Eliahou--Kervaire \cite{E--K} says that,
for a stable ideal $I$ of $A$,
the projective dimension
$\proj \dim A / I$ of the quotient algebra $A / I$ 
coincides with
$\max \{ m(u) : u \in G(I) \}$.
Since a lexsegment ideal is stable,
it follows from Proposition \ref{lexsegment} (a)
together with the Auslander--Buchsbaum formula
\cite[Theorem 1.3.3]{BruHer} that

\begin{Corollary}  
\label{depth}
Let $I$ be a lexsegment ideal of $A$ and $\depth A / I$
the depth of the quotient algebra $A / I$ of $A$.
Then $\depth A / I = \max\{ n - |G(I)|, 0 \}$.
\end{Corollary}

It is known 
that, given a numerical function 
$H : \NN \to \NN$ satisfying $H(0) = 1$, $H(1) \leq n$
and $H(q+1) \leq H(q)^{\langle q \rangle}$ 
for $q = 1, 2, \ldots$, there exists the unique 
lexsegment ideal $I$ of $A$ with $H_{A/I} = H$.
We say that
a numerical function 
$H : \NN \to \NN$ satisfying $H(0) = 1$, $H(1) \leq n$
and $H(q+1) \leq H(q)^{\langle q \rangle}$ 
for $q = 1, 2, \ldots$ is {\em critical} if
the lexsegment ideal $I$ of $A$ with $H_{A/I} = H$
is universal lexsegment.

\begin{Proposition}
\label{hilbert}
A numerical function 
$H : \NN \to \NN$
satisfying $H(0) = 1$, $H(1) \leq n$
and $H(q+1) \leq H(q)^{\langle q \rangle}$ 
for $q = 1, 2, \ldots$ is critical  
if and only if
there is an integer $1 \leq \delta \leq n$ together with
a sequence of integers $(e_1, \ldots, e_\delta)$ with
$1 \leq e_1 \leq \cdots \leq e_\delta$ such that
\begin{eqnarray}
\label{Boston}
H(q) = {n - 1 + q \choose n - 1} - 
\sum_{i=1}^{\delta} {n - i + q - e_i \choose n - i} 
\end{eqnarray}
for $q = 0, 1, \ldots$.
Moreover, $\delta$ is equal to the number of minimal monomial generators of the 
universal lexsegment ideal $I$ of $A$ with $H_{A/I}=H$.
\end{Proposition}

\begin{proof}
First, to prove the ``only if'' part, 
let $I$ be a universal lexsegment ideal of $A$ with
$G(I) = \{ u_1, \ldots, u_\delta \}$, where
$\delta \leq n$.  Suppose that
$\deg u_1 \leq \cdots \leq \deg u_\delta$ 
and that
$u_{i+1} <_{\lex} u_i$
if $\deg u_i = \deg u_{i+1}$.
Proposition \ref{lexsegment} (a) says that,
for $1 \leq i < j \leq \delta$,
the monomial $x_i u_j$ is divided by $u_i$
and no monomial belongs to
both $u_i K[x_i, \ldots, x_n]$ and
$u_j K[x_j, \ldots, x_n]$. 
Hence the direct sum decomposition
$
I = \bigoplus_{i=1}^{\delta} u_i K[x_i, \ldots, x_n]
$
arises.  
Let $e_i = \deg u_i$ for $i = 1, 2, \ldots,\delta$.
The fact that the number of monomials of degree $q$
belonging to $I$ is
$\sum_{i=1}^{\delta} {n - i + q - e_i \choose n - i}$
yields the formula (\ref{Boston}),
as required. 

Next we consider the ``if'' part.
Let $H:\NN \to \NN$ be a numerical function of the form (\ref{Boston}).
Since $1 \leq e_1 \leq \cdots \leq e_\delta$ and $\delta \leq n$,
Proposition \ref{lexsegment} (b) and Corollary \ref{universal}
say that there exists the universal lexsegment ideal
$I$ with $G(I)=\{u_1,\dots,u_\delta\}$ such that $\deg(u_i)=e_i$ for all $i$.
Then the computation of Hilbert functions
in the proof of the ``only if'' part implies $H(I,q)=H(q)$ for all $q \in \NN$.
%
\end{proof}

A {\em critical} ideal of $A$ is a homogeneous ideal
$I$ of $A$ with $I \neq A$ such that the Hilbert function
$H_R$ of the quotient algebra $R = A / I$ is critical.
In other words, a critical ideal of $A$ is a homogeneous 
ideal $I$ of $A$ such that the lexsegment ideal 
$I^{\lex}$ is universal lexsegment, 
where $I^{\lex}$ is the unique lexsegment ideal
of $A$ such that $A / I$ and $A / I^{\lex}$
have the same Hilbert function.
Somewhat surprisingly, 

\begin{Theorem}
\label{critical}
Suppose that a homogeneous ideal $I$ of $A$
is critical.  Then 
\[
\depth A / I = \depth A / I^{\lex}.
\]
\end{Theorem}

\begin{proof}
Let $\beta_{ij}$ (resp.\ $\beta'_{ij}$)
denote the graded Betti numbers of 
$I$ (resp.\ $I^{\lex}$).
Let $G(I^{\lex}) = \{ u_1, \ldots, u_\delta \}$
with $\delta \leq n$, where
$\deg u_1 \leq \cdots \leq \deg u_\delta$ 
and where
$u_{i+1} <_{\lex} u_i$
if $\deg u_i = \deg u_{i+1}$.
Let $e_i = \deg u_i$ for $i = 1, \ldots, \delta$.
Eliahou--Kervaire \cite{E--K} 
together with Proposition \ref{lexsegment} (a) guarantees that
$\beta'_{i, \delta - 1 + e_\delta} = 0$
unless $i = \delta - 1$ and
$\beta'_{\delta - 1, \delta - 1 + e_\delta} = 1$.
Since $A / I$ and $A / I^{\lex}$ 
have the same Hilbert function,
it follows from \cite[Lemma 4.1.13]{BruHer} that
\[
\sum_{i \geq 0}
( - 1 )^i \beta_{i, \delta - 1 + e_\delta}
= \sum_{i \geq 0}
( - 1 )^i \beta'_{i, \delta - 1 + e_\delta}.
\]
Since $\beta_{ij} \leq \beta'_{ij}$
for all $i$ and $j$
(\cite{Big}, \cite{Hul} and \cite{P}),
it follows that 
$\beta_{\delta - 1, \delta - 1 + e_\delta} = 1$.
Thus in particular
$\proj \dim A / I \geq \delta$.
Since 
$\proj \dim A / I^{\lex} = \delta$
and
$\proj \dim A / I \leq \proj \dim A / I^{\lex}$,
it follows that
$\proj \dim A / I = \proj \dim A / I^{\lex}
= \delta$.
Thus we have
$\depth A / I = \dim A / I^{\lex} = n - \delta$,
as desired.
\end{proof}

Moreover,
in case of monomial ideals,
the graded Betti numbers of a critical ideal are determined by its Hilbert function.

\begin{Corollary} \label{omake}
Suppose that a monomial ideal $I$ of $A$ is critical.
Then $I$ and $I^{\lex}$ have the same graded Betti numbers.
\end{Corollary}

\begin{proof}
It follows from Taylor's resolution of monomial ideals (see \cite[p.\ 18]{E--K})
that
$$ \mathrm{proj\ dim} (A/I) \leq |G(I)|.$$
On the other hand,
Corollary \ref{depth} and Theorem \ref{critical}
say that
$$ \mathrm{proj \ dim}(A/I)= \mathrm{proj \ dim}(A/{I^{\lex}})=|G(I^{\lex})|.$$
Since the number of elements in $G(I^{\lex})$ is always larger than that of $G(I)$,
we have $|G(I)| = |G(I^{\lex})|$.
This means 
$\sum_{j \geq 0} \beta_{0j}(I)=\sum_{j \geq 0} \beta_{0j}(I^{\lex})$.
Then it follows from \cite[Theorem 1.3]{C} that
$\beta_{ij}(I)=\beta_{ij}({I^{\lex}})$ for all $i$ and $j$.
\end{proof}

We are not sure that Corollary \ref{omake} holds for an arbitrary critical ideal.

\begin{Example}
{\em
Let $I$ be the monomial ideal
$(x_1x_4, x_3x_4)$ of $K[x_1, x_2, x_3, x_4]$.
Since $I^{\lex} = (x_1^2, x_1x_2)$
is universal lexsegment, it follows that
$\depth A / I = 2$.
}
\end{Example}

\section{Depth and Hilbert functions}
Let, as before, 
$A = K[x_1, \ldots, x_n]$ denote the polynomial ring
in $n$ variables over a field $K$ with each $\deg x_i = 1$.
Given a numerical function 
$H : \NN \to \NN$ satisfying $H(0) = 1$, $H(1) \leq n$
and $H(q+1) \leq H(q)^{\langle q \rangle}$ 
for $q = 1, 2, \ldots$,
we write $\AA_H$ 
for the set of those integers $0 \leq r \leq n$ 
such that there exists a homogeneous ideal $I$ of $A$
with $H_{A/I} = H$ and with
$\depth A / I = r$.
It follows from Corollary \ref{depth}
together with Theorem \ref{critical} that
if $H$ is critical, that is, $H$ is of the form $(\ref{Boston})$,
then $\AA_H = \{ n - \delta \}$.

\begin{Theorem}
\label{noncritical}
Suppose that a numerical function 
$H : \NN \to \NN$  
satisfying $H(0) = 1$, $H(1) \leq n$
and $H(q+1) \leq H(q)^{\langle q \rangle}$ 
for $q = 1, 2, \ldots$ is noncritical.
Then $\AA_H = \{ 0, 1, 2, \ldots, b \}$, where
$b$ is the biggest integer for which $b \in \AA_H$.
\end{Theorem}

\begin{proof}
We may assume that $K$ is infinite.
Let $I$ be a homogeneous ideal of $A$ with $H_{A/I} = H$
and with $\depth A / I = b$.
Let $0 \leq r \leq b$.
Since $K$ is infinite and since $\depth A / I = b$, 
there exists a regular sequence 
$(\theta_1, \ldots, \theta_r)$
of $A / I$ with each $\deg \theta_i = 1$.  
It then follows that there exists a homogeneous ideal $J$ 
of $B = K[x_1, \ldots, x_{n-r}]$ such that
the ideal $J A$ of $A$ 
satisfies $H_{A/(J A)} = H$.

We now claim that the lexsegment ideal $J^{\lex} \subset B$
of $J$ cannot be universal lexsegment.  
In fact, if $J^{\lex}$ is universal lexsegment,
then $J^{\lex}$ remains being lexsegment in
the polynomial ring $K[x_1, \ldots, x_m]$
for each $m \geq n - r$.
In particular the ideal $J^{\lex} A$ of $A$
is universal lexsegment.
Since $H_{A/(JA)} = H_{A/(J^{\lex} A)} = H$,
the numerical function $H$ is critical,
a contradiction.

Since the lexsegment ideal $J^{\lex}$
of $J$ cannot be universal lexsegment,
it follows from Corollaries \ref{universal} and
\ref{depth} that
$\depth B / J^{\lex} = 0$.  
Thus
$\depth A/(J^{\lex} A) = r$.  Hence $r \in \AA_H$,
as desired.
\end{proof}

One may ask a way to compute the integer $b=\max \AA_H$ from $H$.
This integer $b$ can be determined as follows:
Let $H:\mathbb{N} \to \mathbb{N}$ be a numerical function.
The \textit{differential} $\Delta^1(H)$ of $H$ is the numerical function
defined by $\Delta^1(H)(0)=1$ and $\Delta^1(q)=H(q)-H(q-1)$ for $q \geq 1$.
We define \textit{$p$-th differential} $\Delta^p(H)=\Delta^1 (\Delta^{p-1}(H))$ inductively.

\begin{Theorem} \label{added}
Let $H : \NN \to \NN$ be a numerical function 
satisfying $H(0) = 1$, $H(1) \leq n$
and $H(q+1) \leq H(q)^{\langle q \rangle}$ 
for all $q \geq 1$.
Then one has
\begin{eqnarray}
\ \ \ \ \ \  \max A_H= \max\{ p: \Delta^p(H) \mbox{ satisfies } \Delta^p(H)(q+1) \leq \Delta^p(H)(q)^{\langle q \rangle} \mbox{ for } q \geq 1\} \label{LRM}
\end{eqnarray} 
\end{Theorem}

\begin{proof}
If $p$ is an integer which belongs to the righthand side of (\ref{LRM}),
then there exists a homogeneous ideal $J$ of $B=K[x_1,\dots,x_{n-p}]$
such that $H_{B/J}=\Delta^p (H)$.
Recall that if $M$ is a graded $R$-module and $\vartheta_1,\dots,\vartheta_r$ with each $\deg(\vartheta_i)=1$
is a regular sequence of $M$,
then $H_{M/(\vartheta_1,\dots,\vartheta_r) M}=\Delta^p(H_M)$.
Set $M=A/(JA)$.
Then, since $x_n,x_{n-1},\dots,x_{n-p+1}$ is a regular sequence of $A/(JA)$
and  $M/(x_n,\dots,x_{n-p+1})M \cong B/J$,
we have $H_{A/(JA)}=H$ and $\depth(A/(JA)) \geq p$.
This says that the lefthand side of (\ref{LRM}) is larger than or equal to the righthand side of (\ref{LRM}).

On the other hand,
if there exists a homogeneous ideal $I$ of $A$ such that
$H=H_{A/I}$ and $\depth(A/I)=p$,
then, in the same way as Theorem \ref{noncritical}, there exists a homogeneous ideal $J$ of $B=K[x_1,\dots,x_{n-p}]$
such that $H_{A/(JA)}=H$ and $H_{B/J}=\Delta^p(H)$.
Thus the lefthand side of (\ref{LRM}) is smaller than or equal to the righthand side of (\ref{LRM}).
\end{proof}

\begin{Example}
{\em
Let $I$ be the monomial ideal
$(x_1x_4, x_1x_5, x_2x_5, x_3x_5, x_2x_3x_4)$ 
of $A = K[x_1, x_2, x_3, x_4, x_5]$.
Then
$$I^{\lex}=(x_1^2,x_1x_2,x_1x_3,x_1x_4,x_1x_5^2,x_2^3,x_2^2x_3,
x_2^2x_4^2,x_2^2x_4x_5,x_2^2x_5^3,x_2x_3^4,x_2x_3^3x_4^2).$$
Thus $\depth A / I^{\lex} = 0$ by Corollary \ref{depth}.
Since the Hilbert series 
$\sum_{q=0}^{\infty} H_{A/I}(q) \lambda^q$
of $A / I$ is 
$(1 + 2 \lambda - \lambda^2 - \lambda^3)/( 1 - \lambda)^3$,
it follows from \cite[Corollary 4.1.10]{BruHer} that
the Krull dimension of $A/I$ is $3$ and $3 \not\in \AA_H$.
Since $\depth A / I = 2$, 
one has $\AA_H = \{ 0, 1, 2 \}$.
}
\end{Example}

\end{document}